\documentclass[11pt,a4paper]{amsart}

\usepackage[english]{babel}
\usepackage[T1]{fontenc}
\usepackage[utf8]{inputenc}
\usepackage{kpfonts}
\usepackage{fancyhdr}
\usepackage{amsmath,amsfonts,amssymb}

\newtheoremstyle{gtheorem}
     {3pt}        
     {3pt}        
     {\itshape}           
     {\parindent} 
     {\bfseries}   
     {.}          
     {.5em}       
     {}           

\newtheoremstyle{gdefinition}
     {3pt}        
     {3pt}        
     {}           
     {\parindent} 
     {\bfseries}   
     {.}          
     {.5em}       
     {}           

\newtheoremstyle{gremark}
     {3pt}        
     {3pt}        
     {}           
     {\parindent} 
     {\itshape}   
     {\hspace{0.5em}---}          
     {0.5em}       
     {}           

\theoremstyle{gtheorem}
\newtheorem{theo}{Theorem}[section]
\newtheorem{prop}[theo]{Proposition}

\newtheorem*{theo*}{Theorem}

\theoremstyle{gdefinition}

\theoremstyle{gremark}

\renewcommand{\phi}{\varphi}
\renewcommand{\epsilon}{\varepsilon}

\begin{document}

\title{A natural hermitian metric associated with local universal families of compact Ricci-flat K\"ahler manifolds}
\author{Gunnar Þór Magn{\'{u}}sson}
\thanks{I would like to thank my thesis advisor, Jean-Pierre Demailly, for his patience and invaluable guidance.\\This note is to appear in Comptes Rendus - Mathématique.}
\address{Institut Fourier, 100 rue des Maths, BP 74, 38402 St Martin d'Hères}
\email{gunnar.magnusson@fourier.ujf-grenoble.fr}
\maketitle

\begin{abstract}
  Let $\pi : \mathcal X \to S$ be a local universal family of compact Ricci-flat K\"ahler manifolds over a smooth base $S$. The complexified K\"ahler cones of each fiber of the family form a holomorphic fiber bundle $\mathcal K \to S$. We show that there exists a natural hermitian metric $\omega$ on the fiber product $\mathcal X \times_S \mathcal K$. We then discuss the example of elliptic curves in some detail.
\end{abstract}

\section{Introduction}

Let $X_0$ be a compact K\"ahler manifold with zero first Chern class. Examples of such manifolds include complex tori, Calabi-Yau manifolds, and hyperk\"ahler manifolds. Consider a local universal family $\pi : \mathcal X \to S$ of deformations of $X_0$ over a base $S$, which may be taken to be smooth by the Bogomolov-Tian-Todorov theorem \cite{MR915841}. A. Nannicini \cite{MR829406} and G. Schumacher \cite{MR784145} have constructed a Weil-Petersson metric on polarized families of such manifolds, proven that this metric is K\"ahler, and estimated its curvature.

However, not every family of such manifolds can be polarized, as obstructions arise when the Hodge decomposition of the degree 2 cohomology groups of the fibers varies. An example is given by the family $\mathcal X$ of complex tori of dimension $n \geq 2$ parametrized by $n \times n$ complex matrices $s$ whose imaginary part is positive-definite. Thus the methods of Nannicini and Schumacher do not apply to arbitrary families of compact K\"ahler manifolds with zero first Chern class.

We overcome this difficulty by adding the data of complexified K\"ahler classes of each fiber to the family $\pi : \mathcal X \to S$. The choice of a polarization of $\mathcal X$, if one exists, corresponds to a subvariety of this larger space. In this note we outline the construction of a natural hermitian metric on this space. We then review the case of elliptic curves.

\section{Construction of the hermitian form}

Let $X$ be a compact K\"ahler manifold. A \em complexified K\"ahler class \em on $X$ is a $(1,1)$-cohomology class $\alpha$ whose imaginary part contains a K\"ahler metric. The set of complexified K\"ahler classes on $X$ is denoted by $K^{\mathbb C}(X)$. It is an open cone in $H^{1,1}(X,\mathbb C)$.

Let $\pi : \mathcal X \to S$ be a local universal family of K\"ahler manifolds over a smooth base $S$. Under our hypotheses, there exists a holomorphic vector bundle $p : E^{1,1} \to S$ whose fiber over $s$ is $E^{1,1}_s = H^{1,1}(X_s,\mathbb C)$. It is equipped with a connection $\nabla$.\footnote{If $E^2 \to S$ is the bundle whose fibers are $E^2_s = H^2(X_s,\mathbb C)$, then there is a real-analytic inclusion of $E^{1,1}$ into $E^2$ and a smooth projection $E^2 \to E^{1,1}$. The connection $\nabla$ is induced by these morphisms and the Gauss-Manin connection on $E^2$.}
 We denote by $\mathcal K$ the subset of $E$ consisting of the complexified K\"ahler classes on each fiber. By a theorem of Kodaira-Spencer \cite{MR0115189} the subset $\mathcal K$ is open in $E$. Thus it has the structure of a holomorphic fibration over $S$, and inherets the connection $\nabla$ on $E^{1,1}$.

We work over the fiber product $\mathcal X \times_S \mathcal K \stackrel{\nu}{\longrightarrow} S$, where $\nu$ is induced by $\pi$ and the projection onto $\mathcal X$. As $S$ is smooth the fiber product is a smooth complex manifold. It is also a holomorphic fibration over $S$, and its fiber over $s$ is the product $X_s \times K^{\mathbb C}(X_s)$. A point of $\mathcal X \times_S \mathcal K$ is a triple $(z,\alpha,s)$, where $z$ is a point on the manifold $X_s$, $\alpha$ is a complexified K\"ahler class on $X_s$, and $s$ is a point of $S$.

Denote the projections from $\mathcal X \times_S \mathcal K$ to $\mathcal X$ and $\mathcal K$ by $pr_{\mathcal X}$ and $pr_{\mathcal K}$, respectively. Note that the relative tangent bundle of $\mathcal K$ is $p^* E^{1,1}$ and that $\nu = p \circ pr_{\mathcal K}$. The relative tangent bundle of a fiber product splits holomorphically, so there is a short exact sequence
\begin{equation}
  \label{eq:tangents}
  0 \longrightarrow pr_{\mathcal X}^*T_{\mathcal X/S}
  \oplus \nu^* E^{1,1}
    \longrightarrow T_{\mathcal X \times_S \mathcal K} 
    \stackrel{\nu_*}{\longrightarrow}
    \nu^* T_S \longrightarrow 0
\end{equation}
of holomorphic vector bundles over $\mathcal X \times_S \mathcal K$.

\textbf{Observation 1.} There are natural hermitian metrics on on the left hand and right hand terms of the exact sequence. We sketch their construction:

\textit{Proof: } Restrict the short exact sequence \eqref{eq:tangents} to a point $(z,\alpha,s)$ of the manifold. By Yau's solution of the Calabi conjecture \cite{MR480350} there exists a unique Ricci-flat K\"ahler metric $g$ in $\mathop{\text{Im}} \alpha$. This metric gives a hermitian inner product $g(z)$ on $pr_{\mathcal X}^*T_{\mathcal X/S,(z,\alpha,s)} = T_{X_s,z}$.

The K\"ahler metric $g$ also induces a $L^2$ hermitian inner product on $\nu^*E^{1,1}_{(z,\alpha,s)} = E^{1,1}_s$ via passage to $g$-harmonic forms and integration over $X_s$.

Thirdly, the metric $g$ gives a hermitian inner prodcut on the vector space $H^1(X_s,T_{X_s})$ by the same process as above. This gives rise to a hermitian form on $\nu^*T_{S,(z,\alpha,s)} = T_{S,s}$ via pullback by the Kodaira-Spencer morphism of the family $\pi : \mathcal X \to S$. This form is a metric as the family considered is effective.

Finally note that $g$ is the solution of a non-linear elliptic Monge-Amp\`ere equation which depends smoothly on $\alpha$ and $s$. Thus $g$ varies smoothly with these parameters. As the inner products above vary smoothly with $g$, they define smooth hermitian metrics on each vector bundle. \hfill$\square$

\textbf{Observation 2.} Siu's method of constructing canonical lifts defines a smooth lifting of $\nu^*T_S$ into $pr_{\mathcal X}^*T_{\mathcal X}$ at the given point. A sketch of the construction is as follows:

\textit{Proof: } Let $(z,\alpha,s)$ be a point of $\mathcal X \times_S \mathcal K$. Let $g$ be the Ricci-flat K\"ahler metric in $\mathcal{\text{Im}} \alpha$, and let $\rho$ be the Kodaira-Spencer morphism of the family $\pi : \mathcal X \to S$. We lift a section $\xi$ of $\nu^*T_S$ to a smooth section $\eta$ of $pr_{\mathcal X}^*T_{\mathcal X}$ which satisfies that $\bar \partial \eta|_{X_s}$ is $g$-harmonic. A rigorous construction of this lift can be carried out via the theory of harmonic forms on complex manifolds and the use of a relative Green operator. The details will be in the author's PhD thesis. The section $\eta$ is seen to depend smoothly on $g$ and the parameters $(z,\alpha,s)$, so we obtain a smooth lift of $\nu^*T_S$ into $pr_{\mathcal X}^*T_{\mathcal X}$. \hfill$\square$

Now recall that we have a connection $\nabla$ on $\mathcal K$. This connection induces a smooth lift of $\nu^*T_S$ into $pr_{\mathcal K}^*T_{\mathcal K}$ by standard constructions in differential geometry. Together the two lifts define a smooth lift of $\nu^*T_S$ into $T_{\mathcal X \times_S \mathcal K}$, and thus give a smooth splitting of the short exact sequence \eqref{eq:tangents}. Consequently:

\begin{theo}
  There is a natural smooth hermitian metric $\omega$ on the fiber product $\mathcal X \times_S \mathcal K$.
\end{theo}

We recall that the Ricci-flat metrics in each K\"ahler class are unique. One can use this uniqueness to prove that the metric $\omega$ is compatible with base change $B \to S$, and thus descends to any moduli space.

\section{A simple example}

\subsection{Calculation of the natural metric}

Let us recall some conventions and notations.

An elliptic curve $X$ is a compact complex curve of genus 1. Its canonical bundle is trivial, so it is a compact K\"ahler manifold with zero first Chern class. Choosing a point of $X$ makes the curve into a commutative Lie group. A \em marking \em of an elliptic curve is the choice of a basis of $H_1(X,\mathbb Z)$, or equivalently, the choice of a lattice $\Lambda$ such that $X = \mathbb C / \Lambda$. One can always choose a marking such that $\Lambda = \mathbb Z \oplus s \mathbb Z$ for some $s$ in the upper half-plane $\mathbb H$.

Elliptic curves with the choice of a point are parametrized by the orbifold $\mathcal M_{1,1}$, and marked elliptic curves are parametrized by the half-plane $\mathbb H$. There is a morphism $\mathbb H \to \mathcal M_{1,1}$ given by passage to the quotient by the action of $SL_2(\mathbb Z)$ on $\mathbb H$. The space $\mathbb H$ is a fine moduli space for marked elliptic curves.

There is a universal family $\pi : \mathcal X \to \mathbb H$. It can be constructed as the quotient of the space $\mathbb C \times \mathbb H$ by the action of the group $G = \{ g_{n,m} \,|\, g_{n,m}(z,s) = (z + n + ms, s) \}$. A universal family over $\mathcal M_{1,1}$ is then obtained by quotienting again by $SL_2(\mathbb Z)$, whose action lifts to the family over $\mathbb H$.

The fibration of complexified K\"ahler cones $\mathcal K$ over $\mathbb H$ is trivial. We will identify it with the product $\mathbb H \times \mathbb H$, where the couple $(\alpha, s)$ corresponds to the cohomology class on $X_s$ defined by $\frac{\alpha}{\mathop{\text{Im}} s} \, \frac{i}{2} dz \wedge d\bar z$. The action of the group $G$ on $\mathbb C \times \mathbb H$ extends to the fibration of K\"ahler cones. The extension acts trivially on the K\"ahler cones, so the fiber product $\mathcal X \times_{\mathbb H} \mathcal K$ is the quotient of the space $\mathbb C \times \mathbb H \times \mathbb H$ by the group $\tilde G = \{ g_{n,m} \,|\, g_{n,m}(z,\alpha,s) = (z + n + ms, \alpha, s) \}$. Let us denote the morphism of passage to the quotient by $q$. Calculations show that:

\begin{prop}
  Let $\omega$ be the natural hermitian metric on $\mathcal X \times_{\mathbb H} \mathcal K$. Its pullback to $\mathbb C \times \mathbb H \times \mathbb H$ is
\begin{equation}
  \label{eq:pullback}
  q^* \omega =
  \begin{pmatrix}
    g_{\mathcal X/S} & 0 & g_{\mathcal X/S} \, \overline a \\
    0 & g_{L^2} & 0 \\
    a \, g_{\mathcal X/S} & 0 & 
    g_{WP} + a \, g_{\mathcal X/S} \, \overline a
  \end{pmatrix}
\end{equation}
where $g_{\mathcal X/S}(z,\alpha,s) = \frac{\alpha - \bar \alpha}{s - \bar s}$, $g_{L^2}(z,\alpha,s) = \frac{i}{2} \frac{1}{\alpha - \bar \alpha}$, $g_{WP}(z,\alpha,s) = \frac{i}{2} \frac{\alpha - \bar \alpha}{(s - \bar s)^2}$ and $a(z,\alpha,s) = -\frac{z - \bar z}{s - \bar s}$.
\end{prop}

Direct calculations show that this metric is merely hermitian, not K\"ahler.

\subsection{Fibrations and metric}

Choose a point $s_0 \in \mathbb H$ and let $\mathcal K_0$ be the fiber of $\mathcal K$ over $s_0$. As the fibration $\mathcal K$ is trivial over $\mathbb H$, then all that follows will not depend on the choice of the point $s_0$.

The fiber product $\mathcal X \times_{\mathbb H} \mathcal K$ is now equipped with two holomorphic fibrations $\nu : \mathcal X \times_{\mathbb H} \mathcal K \to \mathbb H$ and $\eta : \mathcal X \times_{\mathbb H} \mathcal K \to \mathcal K_0$. These fibrations are induced by the projection onto each of the copy of the half-plane $\mathbb H$ in the product $\mathbb C \times \mathbb H \times \mathbb H$.

The choice of a fiber of each fibration has a natural interpretation: choosing a fiber of $\nu$ means fixing a complex structure on the 2-torus $\mathbb R^2 / \mathbb Z^2$ and letting the symplectic structures on that torus vary; while choosing a fiber of $\eta$ corresponds to choosing a polarization of the family $\pi : \mathcal X \to \mathbb H$. It turns out that a normalized version of the natural metric $\omega$ is aware of these fibrations.

\begin{theo}
  Denote the volume of $X_s$ with respect to the complexified K\"ahler class $\alpha$ by $\mathop{\text{Vol}}(X_s,\alpha) = \mathop{\text{Im}} \alpha$. The normalized hermitian metric $\omega_U := \frac{1}{\mathop{\text{Vol}}(X_s,\alpha)} \, \omega$ becomes K\"ahler when restricted to any fiber of either $\nu$ or $\eta$.
\end{theo}

The current proof of this theorem is by a direct calculation on the space $\mathbb C \times \mathbb H \times \mathbb H$. One can show that the action of $SL_2(\mathbb Z)$ lifts to the fiber product, that the metric $\omega$ is invariant under this action, and that the fibrations descend to the quotient. Thus the theorem also holds on the fiber product over $\mathcal M_{1,1}$.

\subsection{An isometry}

Inspired by the interpretation of mirror symmetry in the case of elliptic curves, we have constructed a self-diffeomorphism $\phi$ of the fiber product $\mathcal X \times_{\mathbb H} \mathcal K$. This diffeomorphism is induced by the map $(z,\alpha,s) \mapsto (L(\alpha,s)(z), s, \alpha)$ on the space $\mathbb C \times \mathbb H \times \mathbb H$, where $L(\alpha,s)$ is the unique $\mathbb R$-linear endomorphism of $\mathbb C$ which sends the basis $(1,s)$ to $(1,\alpha)$. This endomorphism induces a diffeomorphism between the elliptic curves $X_s$ and $X_\alpha$. The diffeomorphism $\phi$ exchanges the two holomorphic fibrations $\nu$ and $\eta$ constructed above.

\begin{theo}
  Equip the fiber product $\mathcal X \times_{\mathbb H} \mathcal K$ with the normalized hermitian metric $\omega_U$ considered above. Then $\phi$ is an auto-isometry of $\mathcal X \times_{\mathbb H} \mathcal K$.
\end{theo}

The proof of this theorem is again by direct calculations. As before, the action of $SL_2(\mathbb Z)$ on $\mathbb H$ lifts to an action of the fiber product $\mathcal X \times_{\mathbb H} \mathcal K$, but the diffeomorphism $\phi$ is not invariant under this action. Thus we do not get a version of this theorem over $\mathcal M_{1,1}$.

\bibliographystyle{amsplain}
\bibliography{arxiv}

\end{document}